\documentclass[11pt]{amsart}
\usepackage{amsfonts,amsmath,amssymb}

\def\ot{\otimes}

\newtheorem{definition}{Definition}

\newtheorem{remark}[definition]{Remark}

\newtheorem{theorem}[definition]{Theorem}

\def\bea{\begin{eqnarray*}}
\def\eea{\end{eqnarray*}}

\def\ot{\otimes}

\def\bea{\begin{eqnarray*}}
\def\eea{\end{eqnarray*}}

\def\ot{\otimes}

\begin{document}

\title{Stable  equivalence  of Morita type and Frobenius extensions}
\author{M. Beattie}
\address{Department of Mathematics and Computer Science, Mount Allison University,
Sackville, New Brunswick, Canada E4L 1E6 }
\email{mbeattie@mta.ca}
\author{S. Caenepeel}
\thanks{The work on this note was started while the last two authors were visiting the Mathematics Department at Mount Allison University. They thank the department for its warm hospitality}
\address{Faculteit Ingenieurswetenschappen,
Vrije Universiteit Brussel,
Pleinlaan 2,
B-1050 Brussel, Belgium}
\email{scaenepe@vub.ac.be}
\author{\c{S}. Raianu}
\address{Mathematics Department, California State University, Dominguez Hills,
1000 E Victoria St, Carson, CA 90747, USA}
\email{sraianu@csudh.edu}
\begin{abstract} A.S. Dugas and R. Mart\'{i}nez-Villa proved in \cite[Corollary 5.1]{dm} that if there exists a stable equivalence of Morita type between the $k$-algebras $\Lambda$ and $\Gamma$, then it is possible to replace $\Lambda$ by a Morita equivalent $k$-algebra $\Delta$ such that $\Gamma$ is a subring of $\Delta$ and the induction and restriction functors induce inverse stable equivalences. In this note we give an affirmative answer to a question of Alex Dugas about the existence of a $\Gamma$-coring structure on $\Delta$. We do this by showing that $\Delta$ is a Frobenius extension of $\Gamma$.
\end{abstract}
\maketitle

As in \cite{dm}, we will assume throughout that  the algebras $\Lambda$ and $\Gamma$ are finite dimensional over a field $k$ and have no semisimple blocks.

The algebras $\Lambda$ and $\Gamma$ are said to be stably equivalent if the categories of finitely generated modules modulo projectives for $\Lambda$ and $\Gamma$ are equivalent (see \cite{AR}).

A pair of left-right projective bimodules ${_{\Lambda}M}_{\Gamma}$ and  ${_{\Gamma}N}_{\Lambda}$ is said to induce a stable equivalence of Morita type between $\Lambda$ and $\Gamma$ if we have the following isomorphisms of bimodules:
$${_{\Lambda}M}\otimes_{\Gamma}N_{\Lambda}\simeq{_{\Lambda}\Lambda}_{\Lambda}\oplus {_{\Lambda}P}_{\Lambda}\mbox{\;\;\;{\it and}\;\;\;}{_{\Gamma}N}\otimes_{\Lambda}M_{\Gamma}\simeq{_{\Gamma}\Gamma}_{\Gamma}\oplus{_{\Gamma}Q}_{\Gamma}$$
where   ${_{\Lambda}P}_{\Lambda}$ and ${_{\Gamma}Q}_{\Gamma}$ are projective bimodules (see \cite{b}).

We begin by stating the result of Dugas and Mart\'{i}nez-Villa mentioned in the abstract:
\begin{theorem}\label{cor51} (see \cite[Corollary 5.1]{dm})
Let $\Lambda$ and $\Gamma$ be finite dimensional $k$-algebras whose semisimple quotients are separable. If at least one of them is indecomposable, then the following are equivalent:\\
(1) There exists a stable equivalence of Morita type between $\Lambda$ and $\Gamma$.\\
(2) There exists a $k$-algebra $\Delta$, Morita equivalent to $\Lambda$, and an injective ring homomorphism
$\Gamma\hookrightarrow\Delta$ such that the restriction and induction functors are exact and induce inverse stable equivalences.\\
(3) There exists a $k$-algebra $\Delta$, Morita equivalent to $\Lambda$, and an injective ring homomorphism
$\Gamma\hookrightarrow\Delta$ such that
$${_{\Gamma}\Delta}_{\Gamma}={_{\Gamma}\Gamma}_{\Gamma}\oplus {_{\Gamma}P}_{\Gamma}\mbox{\;\;\;
{\it and}\;\;\;}{_{\Delta}\Delta}\otimes_{\Gamma}\Delta_{\Delta}\simeq{_{\Delta}\Delta}_{\Delta}\oplus{_{\Delta}Q}_{\Delta}$$
for projective bimodules  ${_{\Gamma}P}_{\Gamma}$ and ${_{\Delta}Q}_{\Delta}$.
\end{theorem}

We recall now the definition of Frobenius extension, and its dual notion, Frobenius coring.

\begin{definition}(see \cite{k})  Let $i:R\longrightarrow S$ be a ring homomorphism. Then  $S/R$ is called a Frobenius extension if one of the following equivalent conditions is satisfied:
\begin{enumerate}
\item $S$ is finitely generated and projective as a right $R$-module and $Hom_R(S,R)$ and $S$ are isomorphic as $(R,S)$-bimodules;
\item there exists a Frobenius system $(e,\varepsilon)$, consisting of 
$$e=e^1\ot e^2\in (S\ot_R S)^S=\{e^1\ot e^2\in S\ot_R S\mid se^1\ot e^2=e^1\ot e^2s, \forall s\in S\}$$
and $\varepsilon:\ S\to R$ an $R$-bimodule map such that $\varepsilon(e^1)e^2=e^1\varepsilon(e^2)=1$.
\end{enumerate}
\end{definition}

For the proof of the equivalence of the two conditions, see for example \cite[Theorem 28]{cmz}.

\begin{definition}(see \cite{s}) If $R$ is a ring, a coring is a comonoid in the monoidal category of $R$-bimodules. So a coring consists of an $R$-bimodule $C$, together with a coassociative  comultiplication $C\longrightarrow C\otimes_RC$ and counit $C\longrightarrow R$ which are both $R$-bimodule maps.\\
$C$ is called a Frobenius $R$-coring if there exists a Frobenius system $(\theta, 1)$,
consisting of an element $1\in C$ and an $R$-bimodule map $\theta:\ C\ot_RC\to R$ satisfying the
conditions
$$c_{(1)}\theta(c_{(2)}\ot d)=\theta(c\ot d_{(1)})d_{(2)}~~{\rm and}~~\theta(c\ot 1)=\theta(1\ot c)=
\varepsilon(c).$$
\end{definition}

Let $(S,m,1,e,\varepsilon)$ be a Frobenius extension of $R$, and consider $\Delta:\
S\to S\ot_RS$, $\Delta(s)=se=es$. An easy verification shows that $(S,\Delta,\varepsilon,\theta=
\varepsilon\circ m,1)$ is a Frobenius coring.\\
Conversely, if $(C,\Delta,\varepsilon,\theta,1)$ is a Frobenius $R$-coring, then 
$(C,m,1,\Delta(1),\varepsilon)$ is a Frobenius extension. Here $m:\ C\ot_RC\to C$,
$m(c\ot d)=c_{(1)}\theta(c_{(2)}\ot d)=\theta(c\ot d_{(1)})d_{(2)}$.\\
These two assertions basically tell us that Frobenius extension structures on an $R$-bimodule $M$
correspond bijectively to Frobenius $R$-coring structures on $M$.\\
Let $S$ be a Frobenius extension. Then the categories ${\mathcal M}_S$ and ${\mathcal M}^S$ are
isomorphic: on a right $S$-module, we define a right $S$-coaction by
$\rho(m)=me^1\ot e^2$. On a right $S$-comodule, we define a right $S$-action 
$ms=m_{[0]}\varepsilon(m_{[1]}s)$.\\
The restriction functor $G: {\mathcal M}_S\to {\mathcal M}_R$ has a left adjoint, the induction functor $F$;
the forgetful functor ${\mathcal M}_S\to {\mathcal M}_R$ has a right adjoint. These functors are
compatible with the above isomorphism. This implies that $G$ is at the same time a left and a
right adjoint of $F$.

\begin{definition}\label{frobfunc} (see \cite{m} or \cite[p.91]{cmz})
Let $F:{\mathcal C}\longrightarrow{\mathcal D}$ be a covariant functor. If there exists a functor $G:{\mathcal D}\longrightarrow{\mathcal C}$  which is at the same time a right and a left adjoint of $F$, then we call $F$ a {\em Frobenius functor}, and we say that $(F,G)$ is a {\em Frobenius pair} for ${\mathcal C}$ and ${\mathcal D}$.
\end{definition}

\begin{remark}\label{frobfuncext} (see \cite{k} or \cite[Theorem 28, p.103]{cmz}) Let $i:R\longrightarrow S$ be a ring homomorphism, $F$ the induction functor and $G$ the restriction functor. If $S/R$ is a Frobenius extension,
then we have seen above that $(F,G)$ is a Frobenius pair; in fact, it can be shown that the converse
also holds: $(F,G)$ is a Frobenius pair if and only if $S/R$ is a Frobenius extension.
\end{remark}

We can now state and prove our result. Assertion (3) gives an affirmative answer to a question asked by Alex Dugas.
\begin{theorem}\label{frobext}
Let $\Lambda$ and $\Gamma$ be finite dimensional $k$-algebras whose semisimple quotients are separable. Assume that at least one of them is indecomposable, and that there exists a stable equivalence of Morita type between $\Lambda$ and $\Gamma$. Then the following assertions hold:\\
(1) There exists a $k$-algebra $\Delta$, Morita equivalent to $\Lambda$, and an injective ring homomorphism $\Gamma\hookrightarrow\Delta$ such that the restriction and induction functors are a Frobenius pair.\\
(2) There exists a $k$-algebra $\Delta$, Morita equivalent to $\Lambda$, and an injective ring homomorphism $\Gamma\hookrightarrow\Delta$ such that $\Delta/\Gamma$ is a Frobenius extension.\\
(3) There exists a $k$-algebra $\Delta$, Morita equivalent to $\Lambda$, and an injective ring homomorphism $\Gamma\hookrightarrow\Delta$ such that
$${_{\Gamma}\Delta}_{\Gamma}={_{\Gamma}\Gamma}_{\Gamma}\oplus {_{\Gamma}P}_{\Gamma}\mbox{\;\;\;{\it and}\;\;\;}{_{\Delta}\Delta}\otimes_{\Gamma}\Delta_{\Delta}\simeq{_{\Delta}\Delta}_{\Delta}\oplus{_{\Delta}Q}_{\Delta}$$
for projective bimodules  ${_{\Gamma}P}_{\Gamma}$ and ${_{\Delta}Q}_{\Delta}$, and $\Delta$ is a 
Frobenius $\Gamma$-coring with comultiplication given by the injection of ${_{\Delta}\Delta}_{\Delta}$ into ${_{\Delta}\Delta}\otimes_{\Gamma}\Delta_{\Delta}$, and counit given by the projection of ${_{\Gamma}\Delta}_{\Gamma}$ onto ${_{\Gamma}\Gamma}_{\Gamma}$.\\
\end{theorem}
\begin{proof}
(1) Suppose ${_{\Lambda}M}_{\Gamma}$ and ${_{\Gamma}N}_{\Delta}$ are indecomposable bimodules that induce a stable equivalence of Morita type. Let $\Delta=End_{\Lambda}(M)$. By the proof of (1)$\Rightarrow$(2) of  \cite[Corollary 5.1]{dm}, we have that
$$Res^{\Delta}_{\Gamma}\simeq(-\otimes_{\Lambda}M_{\Gamma})\circ Hom_{\Delta}(M,-)$$
and
 $$Ind^{\Delta}_{\Gamma}\simeq(-\otimes_{\Lambda}M_{\Delta})\circ(-\otimes_{\Gamma}N_{\Lambda}).$$
Now $-\otimes_{\Lambda}M_{\Gamma}$ is a right and left adjoint of $-\otimes_{\Gamma}N_{\Lambda}$ by \cite[Corollary 3.1,(2)]{dm}, and $Hom_{\Delta}(M,-)$ is a right and left adjoint of $-\otimes_{\Lambda}M_{\Delta}$ because they are inverse equivalences, so $Res^{\Delta}_{\Gamma}$ is a right and left adjoint of $Ind^{\Delta}_{\Gamma}$.\\
(2) follows from (1) and Remark \ref{frobfuncext}.\\
(3) follows immediately from the above observation that a Frobenius extension is also a Frobenius coring.
\end{proof}



\thebibliography{MMM}



\bibitem{AR}
M. Auslander, I. Reiten, {\it Stable equivalences of artin algebras}, Proc. Conf. on Orders, Group Rings and Related Topics, Lecture Notes in Math., vol {\bf 353}, Springer-Verlag, New York, 1973.

\bibitem{b}
M. Brou\'{e}, Equivalences of blocks of group algebras, in:{\it Finite Dimensional Algebras and related Topics}, Kluwer, 1994, 1-26.



\bibitem{cmz}
S. Caenepeel, G. Militaru, S. Zhu, {\it Frobenius and separable functors for generalized module categories and nonlinear equations} Lecture Notes in Mathematics {\bf 1787} Springer-Verlag, Berlin, 2002.


\bibitem{dm}
A.S. Dugas, R. Mart\'{i}nez-Villa, A note on stable equivalences of Morita type, J. Pure Appl. Algebra {\bf 208} (2007), no. 2, 421-433.





\bibitem{k}
F. Kasch, Projektive Frobenius-Erweiterungen, S.-B. Heidelberger Akad. Wiss.(Math.-Nat. Kl.) {\bf 61} (1961), 89-109.


\bibitem{m}
K. Morita, Adjoint pairs of functors and Frobenius extensions, Sci. Rep. Tokyo Kyoiku Daigaku (Sect. A) {\bf 9} (1965), 40-71.





\bibitem{s}
M. Sweedler, The predual theorem to the Jacobson-Bourbaki theorem,
Trans. Amer. Math. Soc. {\bf 213} (1975), 391–406.
\end{document}